\documentclass[11pt]{article}
\usepackage{latexsym}
\usepackage{epsfig,amssymb,amsfonts,amsmath,amsthm}

\newcommand{\Ker}{\mathcal{K}}
\newcommand{\bw}{\mathbf{w}}
\newcommand{\La}{\mathcal{L}}

\newcommand{\eps}{\varepsilon}
\newcommand{\Pro}[1]{\mathbb{P} \left[\,#1\,\right]}
\newcommand{\Ex}[1]{\mathbb{E} \left[\, #1\,\right]}
\newcommand{\Be}{\mathsf{Be}}
\newcommand{\Bin}{\mathsf{Bin}}
\newcommand{\Var}{\mathsf{Var}}

\newcommand{\remove}[1]{}

\renewcommand{\leq}{\leqslant}
\renewcommand{\geq}{\geqslant}
\renewcommand{\le}{\leqslant}
\renewcommand{\ge}{\geqslant}

\renewcommand{\epsilon}{\varepsilon}

\renewcommand{\bw}{\mathbf{w}}

\newcommand{\Ao}{\mathcal{A}_0}
\newcommand{\Af}{\mathcal{A}_f}

\newtheorem{theorem}{Theorem}  

\newtheorem{lemma}[theorem]{Lemma}

\newtheorem{claim}[theorem]{Claim}
\newtheorem{corollary}[theorem]{Corollary}

\newtheorem{proposition}[theorem]{Proposition}

\newtheorem{definition}[theorem]{Definition}

\numberwithin{theorem}{section}
\numberwithin{equation}{section}

\title{Bootstrap percolation in power-law random graphs}

\author{ Hamed Amini\footnote{Financial support from the Austrian Science Fund (FWF) project P21709 is gratefully acknowledged.} \\ \small{\'Ecole Polytechnique F\'ed\'erale de Lausanne} \\ \small{1015 Lausanne, Switzerland} \\ \small{Hamed.Amini@epfl.ch}\and Nikolaos Fountoulakis\footnote{This
research has been supported by a Marie Curie Intra-European Research Fellowship PIEF-GA-2009-255115 hosted by the Max-Planck-Institut f\"ur
Informatik.} \\
\small{School of Mathematics} \\ \small{University of Birmingham} \\ \small{Edgbaston B15 2TT, UK} \\ \small{N.Fountoulakis@bham.ac.uk}}
\date{}

\begin{document}

\maketitle

\begin{abstract}
A bootstrap percolation process on a graph $G$ is an ``infection" process which evolves in rounds. Initially, there is a subset
of infected nodes and in each subsequent round each uninfected node which has at least $r$ infected neighbours becomes infected
and remains so forever. The parameter $r\geq 2$ is fixed. Such processes have been used as models for the spread of ideas or trends
within a network of individuals.

We analyse this process in the case where the underlying graph is an inhomogeneous random graph, which exhibits a power-law degree
distribution, and initially there are $a(n)$ randomly infected nodes. The main focus of this paper is the number of vertices that will
have been infected by the end of the process. The main result of this work is that if the degree sequence of the random graph follows a power
law with exponent $\beta$, where $2 < \beta < 3$, then a sublinear  number of initially infected vertices is enough to spread the infection
over a linear fraction of the nodes of the random graph, with high probability.

More specifically, we determine explicitly a critical function $a_c(n)$ such that $a_c(n)=o(n)$ with the following property. Assuming that $n$ is
the number of vertices of the underlying random graph, if $a(n) \ll a_c(n)$, then the process does not evolve at all, with high probability as
$n$ grows, whereas if $a(n)\gg a_c(n)$, then there is a constant $\eps>0$ such that, with high probability, the final set of infected vertices
has size at least $\eps n$. This behaviour is in sharp  contrast with the case where the underlying graph is a $G(n,p)$ random graph
with $p=d/n$.
Recent results of Janson, {\L}uczak, Turova and Vallier have shown that if the number of initially infected vertices is sublinear, then
there is lack of evolution of the process.

It turns out that when the maximum degree is $o(n^{1/(\beta -1)})$, then $a_c(n)$ depends also on $r$. But
when the maximum degree is $\Theta (n^{1/(\beta -1)})$, then $a_c (n)=n^{\beta -2 \over \beta -1}$.

\vspace{0.4cm}

\vspace{0.2cm}

\end{abstract}
\newpage
\section{Introduction}

Bootstrap percolation was introduced by Chalupa, Leath and Reich~\cite{ChLeRe:79} in 1979 in the context of
magnetic disordered systems and has been re-discovered since then by several authors mainly due to its connections
with various physical models. A \emph{bootstrap percolation process} with \emph{activation threshold} an integer $r\geq 2$
on a graph $G=G(V,E)$ is a deterministic process which evolves in rounds.
Every vertex has two states: it is either \emph{infected} or \emph{uninfected}.
Initially, there is a subset $\Ao \subseteq V$ which consists of infected vertices, whereas every other vertex is uninfected.
This set can be selected either deterministically or randomly.
Subsequently, in each round, if an uninfected vertex has at least $r$ of its neighbours infected, then it also becomes infected
and remains so forever. This is repeated until no more vertices become infected. We denote the final infected set by $\Af$.

Bootstrap percolation processes (and extensions) have been  used as models
to describe several complex phenomena in diverse areas, from jamming transitions~\cite{tobifi06} and magnetic systems~\cite{sadhsh02} to
neuronal activity~\cite{Am-nn, ET09}. 
Bootstrap percolation also has connections to the dynamics of the Ising model at zero
temperature~\cite{Fontes02}. A short survey regarding applications of bootstrap percolation processes can be found in~\cite{AdL03}.

These processes have also been studied on a variety of graphs, such as trees~\cite{BPP06, FS08}, grids~\cite{CM02, holroyd03, BBDM2010},
hypercubes~\cite{BB06}, as well as on several distributions of random graphs~\cite{balpit07, ar:JLTV10, Am-bp}. In particular, consider the
case when $G$ is the two-dimensional grid on $[n]^2=\{1, \dots, n\}^2$ (i.e., a finite square $[n]^2$ in the square lattice), and $r=2$ (i.e., an uninfected site becomes infected if at least
two of its four neighbours are infected). Then, for an initial set $\Ao \subseteq V$ whose elements are chosen independently at random, each
with probability $p(n)$, the following sharp threshold was determined by Holroyd \cite{holroyd03}. The probability $I(n, p)$ that the entire
square is eventually infected satisfies $I(n,p) \rightarrow 1$ if $\liminf_{n \rightarrow \infty} p(n) \log n > \pi^2/18$,
and  $I(n,p) \rightarrow 0$ if $\limsup_{n \rightarrow \infty} p(n) \log n < \pi^2/18$.
A generalization of this result to the higher dimensional case has been recently proved by Balogh, Bollob\`as and Morris \cite{BBM09}
(when $G$ is the 3-dimensional grid on $[n]^3$ and $r=3$) and Balogh, Bollob\`as, Duminil-Copin and Morris \cite{BBDM2010} (in general).

In the context of real-world networks and in particular in social networks, a bootstrap percolation process can be thought of as a primitive
model for the spread of ideas or new trends within a set of individuals which form a network. Each of them has a threshold $r$ and
$\Ao$ corresponds to the set of individuals who initially are ``infected" with a new belief. If for an ``uninfected" individual at least
$r$ of its acquaintances have adopted the new belief, then this individual adopts it as well.

More than a decade ago, Faloutsos et al.\ \cite{ar:fff99} observed that the Internet exhibits a \emph{power-law} degree distribution,
meaning that the proportion of vertices of degree~$k$ scales like $k^{-\beta}$, for all
sufficiently large~$k$, and some~$\beta > 2$.
In particular, the work of Faloutsos et al.~\cite{ar:fff99} suggested that the degree distribution of the Internet at the
router level follows a power law with $\beta \approx 2.6$. Kumar et al.~\cite{Kum} also provided evidence on the degree distribution of the
World Wide Web viewed as a directed graph on the set of web pages, where a web page ``points" to another web page if the former contains
a link to the latter. They found that the indegree distribution follows a power law with exponent approximately 2.1, whereas the outdegree
distribution follows also a power law with exponent close to 2.7. Other empirical evidence on real-world networks has provided examples
of power law degree distributions with exponents between 2 and 3, see e.g.,~\cite{jtaob00, AlBa2002}.

Thus, in the present work, we focus on the case where $2 < \beta  <3$. More specifically, the underlying random
graph distribution we consider was introduced by Chung and Lu~\cite{CL02},
who invented it as a general purpose model for generating graphs with a
power-law degree sequence. Consider the vertex set $[n]:=\{1,\ldots, n\}$.
Every vertex~$i \in [n]$ is assigned a positive weight~$w_i$, and the pair~$\{i, j\}$, for $i\not = j \in [n]$,
is included in the graph as an edge with probability proportional to~$w_i w_j$, independently of every other pair.
Note that the expected degree of~$i$ is close to~$w_i$. With high probability the degree sequence of the resulting graph follows
a power law, provided that the sequence of weights follows a power law (see~\cite{book:vdH} for a detailed discussion).
Such random graphs are also characterized as \emph{ultra-small worlds}, due to the fact that the typical distance of two vertices that
belong to the same component is $O(\log \log n)$ -- see~\cite{CL03} or~\cite{book:vdH}.

Regarding the initial conditions of the bootstrap percolation process,
our general assumption will be that the initial set of infected vertices $\Ao$ is chosen randomly among
all subsets of vertices of a certain size.

The aim of this paper is to analyse the evolution of the bootstrap percolation process on such random graphs and, in particular,
the typical value of the ratio $|\Af|/ |\Ao|$.
The main finding of the present work is the existence of a critical function $a_c (n)$ such that when $|\Ao|$ ``crosses" $a_c(n)$ we have
a sharp change on the evolution of the bootstrap percolation process. When $|\Ao| \ll a_c (n)$, then typically the process does not evolve, but
when $|\Ao |\gg a_c(n)$, then a linear fraction of vertices is eventually infected. Note that $|\Ao|$ itself may be sublinear.
What turns out to be the key to such a dissemination of the infection is the vertices of high weight. These are typically the vertices that
have high degree in the random graph and, moreover, they form a fairly dense graph. We exploit this fact and show how this causes the
spread of the infection to a linear fraction of the vertices (see Theorem~\ref{thm:fkernel} below).
Interpreting this from the point of view of a social network, these vertices correspond to popular and attractive individuals with many
connections -- these are the \emph{hubs} of the network.
Our analysis sheds light to the role of these individuals in the infection process.

These results are in sharp contrast with the behaviour of the bootstrap percolation process in $G(n,p)$ random graphs, where every edge on a set of $n$ vertices is included independently with probability $p$.
Recently, Janson, {\L}uczak, Turova and Vallier~\cite{ar:JLTV10} came up with a complete analysis of the bootstrap percolation process for various ranges
of the probability $p$. Since the random graphs we consider have constant
average degree, we focus on their findings regarding the range where $p = d/n$ and $d>0$ is fixed.
Among the findings of Janson et al.~\cite{ar:JLTV10} (see Theorem 3.1(i)) is that when $|\Ao| = o(n)$, then typically no evolution occurs.
In other words, the density of the initially infected vertices  must be positive in order for the density of infected vertices to grow. We note that
similar behavior to the case of $G(n,p)$ has been observed in the case of random regular graphs \cite{balpit07}, and in random graphs with
given vertex degrees constructed by configuration model, studied by the first author in \cite{Am-bp}, when the sum of the square of degrees
scales linearly with $n$, the size of the graph. The later case corresponds more to random graphs with power-law degree sequence with
exponent $\beta > 3$. Our results show that the two regimes  $2 < \beta < 3$ and  $\beta > 3$ have completely different behaviors.

\paragraph{Basic notation.}
Let $\mathbb{R}^+$ be the set of positive real numbers.
For non-negative sequences $x_n$ and $y_n$, we describe their relative order of magnitude using Landau's $o(.)$ and $O(.)$ notation. We
write $x_n = O(y_n)$ if there exist $N \in \mathbb{N}$ and $C > 0$ such that $x_n \leq C y_n$ for all $n \geq N$, and
$x_n = o(y_n)$, if $x_n/ y_n \rightarrow 0$, as $n \rightarrow \infty$.
We also write $x_n \ll y_n$ when $x_n = o(y_n)$.

Let $\{ X_n \}_{n \in \mathbb{N}}$ be a sequence of real-valued random variables on a sequence of probability spaces
$\{ (\Omega_n, \mathbb{P}_n)\}_{n \in \mathbb{N}}$.
If $c \in \mathbb{R}$ is a constant, we write $X_n \stackrel{p}{\rightarrow} c$ to denote that $X_n$ \emph{converges in probability to $c$}.
That is, for any $\eps >0$, we have $\mathbb{P}_n (|X_n - c|>\eps) \rightarrow 0$ as $n \rightarrow \infty$.
\\
Let $\{ a_n \}_{n \in \mathbb{N}}$ be a sequence of real numbers that tends to infinity as $n \rightarrow \infty$.
We write $X_n = o_p (a_n)$, if $|X_n|/a_n$ \emph{converges to 0 in probability}.
Additionally, we write $X_n = O_p (a_n)$, to denote that for any positive-valued function $\omega (n) \rightarrow \infty$,
as $n \rightarrow \infty$, we have $\mathbb{P} (|X_n|/a_n \geq \omega (n)) = o(1)$.
If $\mathcal{E}_n$ is a measurable subset of $\Omega_n$, for any $n \in \mathbb{N}$, we say that the sequence
$\{ \mathcal{E}_n \}_{n \in \mathbb{N}}$ occurs \emph{asymptotically almost surely (a.a.s.)} if $\mathbb{P} (\mathcal{E}_n) = 1-o(1)$, as
$n\rightarrow \infty$.
\\
Also, we denote by $\Be (p)$ a Bernoulli distributed random variable whose probability of being equal to 1 is $p$.
The notation $\Bin (k,p)$ denotes a binomially distributed random variable corresponding to the number of
successes of a sequence of $k$ independent Bernoulli trials each having probability of success equal to $p$.

\section{Models and results}\label{sec:models}


The random graph model that we consider is asymptotically equivalent to a model
considered by Chung and Lu~\cite{CL03}, and is a special case of the so-called
\emph{inhomogeneous random graph}, which was introduced S\"oderberg~\cite{ar:s02} and
was studied in great detail by Bollob\'as, Janson and Riordan in~\cite{ar:bjr07}.

\subsection{Inhomogeneous random graphs -- The Chung-Lu model}

In order to define the model we consider for any $n\in\mathbb{N}$ the vertex
set~$[n] := \{1, \dots, n\}$. Each vertex~$i$ is assigned a positive
weight~$w_i(n)$, and we will write~$\bw=\mathbf{w}(n) = (w_1(n), \dots, w_n(n))$. We
assume in the remainder that the weights are deterministic, and we will suppress the
dependence on $n$, whenever this is obvious from the context. However, note that the weights could be themselves
random variables; we will not treat this case here, although it is very likely
that under suitable technical assumptions our results generalize to this case as
well. For any~$S \subseteq [n]$, set
\[
	W_S(\mathbf{w}) := \sum_{i\in S} w_i.
\]
In our random graph model, the event of including the edge~$\{i,j\}$ in the
resulting graph is independent of the events of including all other edges, and
equals
\begin{equation}
\label{eq:pijCL}
	p_{ij}(\mathbf{w}) = \min\left\{\frac{w_iw_j}{W_{[n]}(\mathbf{w})},
1\right\}.
\end{equation}
This model was considered by Chung et \ al., for fairly general choices of~$\mathbf{w}$,
who studied in a series of papers~\cite{CL02,CL03,ar:cl04} several typical properties of the resulting graphs,
such as the average path length or the component
distribution. We will refer to this model as the \emph{Chung-Lu} model, and we shall
write~$CL(\mathbf{w})$ for a random graph in which each possible edge~$\{i,j\}$ is included
independently with probability as in~\eqref{eq:pijCL}. Moreover, we will suppress the dependence on
$\mathbf{w}$, if it is clear from the context which sequence of weights we refer to.

Note that in a Chung-Lu random graph, the weights essentially control the
\emph{expected} degrees of the vertices. Indeed, if we ignore the minimization
in~\eqref{eq:pijCL}, and also allow a loop at vertex~$i$, then the expected
degree of that vertex is~$\sum_{j=1}^n w_iw_j/W_{[n]} = w_i$. In the general
case, a similar asymptotic statement is true, unless the weights fluctuate too
much. Consequently, the choice of~$\mathbf{w}$ has a significant effect on the
degree sequence of the resulting graph. For example, the authors of~\cite{CL03}
choose~$w_i = d {\beta -2 \over \beta -1}(\frac{n}{i + i_0})^{1/(\beta - 1)}$, which typically results in
a graph with a power-law degree sequence with exponent~$\beta$, average
degree~$d$, and maximum degree proportional
to~$({n}/{i_0})^{1/(\beta - 1)}$, where~$i_0$ was chosen such that this
expression is~$O(n^{1/2})$.
Our results will hold in a more general setting,
where larger fluctuations around a ``strict'' power law are allowed, and also
larger maximum degrees are possible, thus allowing a greater flexibility in
the choice of the parameters.

\subsection{Power-law degree distributions}

Following van der Hofstad~\cite{book:vdH}, let us write for any $n\in \mathbb{N}$ and
any sequence of weights $\mathbf{w} = (w_1(n), \dots, w_n(n))$
\[
	 F_n(x) = n^{-1} \sum_{i=1}^n \mathbf{1}[w_i(n) < x], \ \ \forall x\in [0,\infty)
\]
for the empirical distribution function of the weight of a vertex chosen uniformly at random.
We will assume that $F_n$ satisfies the following two  conditions.
\begin{definition}
\label{def:regcond}
We say that $(F_n)_{n \ge 1}$ is \emph{regular}, if it has the following two
properties.
\begin{itemize}
  \setlength{\itemsep}{1pt}
  \setlength{\parskip}{0pt}
  \setlength{\parsep}{0pt}
  \item {\bf[Weak convergence of weight]} There is a distribution function
$F:[0,\infty)\to [0,1]$ such that for all $x$ at which $F$ is continuous
$\lim_{n\to\infty}F_n(x) = F(x)$;
\item {\bf[Convergence of average weight]} Let $W_n$ be a random variable with
distribution function $F_n$, and let $W_F$ be a random variable with
distribution function $F$. Then we have $\lim_{n\to\infty}\Ex{W_n} = \Ex{W_F}$.
\end{itemize}
\end{definition}
The regularity of~$(F_n)_{n\ge 1}$ guarantees two important properties. Firstly,
the weight of a random vertex is approximately distributed as a random variable that follows a certain distribution.
Secondly, this variable has finite mean and therefore the resulting graph has bounded average degree.
Apart from regularity, our focus will be on weight sequences that give rise to power-law degree distributions.
\begin{definition}
\label{def:Distr}
We say that a regular sequence~$(F_n)_{n \ge 1}$ is \emph{of power law with exponent~$\beta$},
if there are~$0<\gamma_1<\gamma_2$, $x_0 >0$ and $0 < \zeta \leq {1/(\beta -1)} $ such that for all~$x_0 \le x \le n^{\zeta}$
\[
	\gamma_1 x^{-\beta + 1}\le 1 - F_n(x) \le \gamma_2 x^{-\beta + 1},
\]
and $F_n(x) = 0$ for $x < x_0$, but $F_n(x) = 1$ for $x > n^{\zeta}$.
\end{definition}
Thus, we may assume that for $1\leq i \leq n(1-F_n(n^{\zeta}))$ we have $w_i = n^{\zeta}$, whereas for
$(1-F_n(n^{\zeta}))n < i \leq n$ we have $w_i = [1-F_n]^{-1} (i/n)$, where
$[1-F_n]^{-1}$ is the generalized inverse of $1-F_n$, that is,
for $x \in [0,1]$ we define $[1-F_n]^{-1}(x) = \inf \{ s \ : \ 1-F_n(s) < x \}$.
Note that according to the above definition, for $\zeta > {1/(\beta -1)}$, we have $n(1-F_n(n^\zeta)) = 0$, since $1-F_n(n^\zeta) \leq \gamma_2
n^{-\zeta(\beta-1)}= o(n^{-1})$. So it is natural to assume that $\zeta  \leq {1/(\beta -1)}$.
Recall finally that in the Chung-Lu model~\cite{CL03} the maximum weight is $O(n^{1/2})$.

\subsection{Results}
The main theorem of this paper regards random infection of the whole of $[n]$. We determine explicitly a critical function which we
denote by $a_c(n)$ such that when we infect randomly $a(n)$ vertices in $[n]$, then the following threshold phenomenon occurs.
If $a(n) \ll a_c(n)$, then a.a.s. the infection spreads no further than $\Ao$, but when $a(n) \gg a_c(n)$, then at least $\eps n$ vertices
become eventually infected, for some $\eps >0$. We remark that $a_c(n) = o(n)$.

\begin{theorem} \label{thm:global}
For any $\beta \in (2,3)$ and any integer $r\geq 2$, we let $a_c (n) = n^{r(1-\zeta) + \zeta (\beta -1)-1\over r}$ for all $n \in \mathbb{N}$.
Let $a: \mathbb{N} \rightarrow \mathbb{N}$ be a function such that $a(n) \rightarrow \infty$, as $n \rightarrow \infty$, but $a(n)=o(n)$.
Let also ${r-1 \over 2r - \beta + 1} < \zeta \leq {1\over \beta -1}$.
If we initially infect randomly $a(n)$ vertices in $[n]$, then the following holds:
\begin{itemize}
\item if $a(n) \ll a_c(n)$, then a.a.s. $\Af = \Ao$;
\item if $a(n) \gg a_c (n)$, then there exists $\eps >0$ such that a.a.s. $|\Af| > \eps n$.
\end{itemize}
When $0< \zeta \leq {r-1 \over 2r - \beta + 1}$ setting $ a_c^+ (n) =n^{1- \zeta~{r-\beta +2 \over r-1 }}$, the following holds.
\begin{itemize}
\item if $a(n) \ll a_c(n)$, then a.a.s. $\Af = \Ao$;
\item if $a(n) \gg a_c^+ (n)$, then there exists $\eps >0$ such that a.a.s. $|\Af| > \eps n$.
\end{itemize}
\end{theorem}

\medskip

Note that the above theorem implies that when the maximum weight of the sequence is $n^{1/(\beta -1)}$, then
the threshold function becomes equal to $n^{\beta - 2 \over \beta -1}$ and does not depend on $r$.

\medskip

The second theorem has to do with the targeted infection of  $a(n)$ vertices where $a(n) \rightarrow \infty$, as $n \rightarrow \infty$.
Let $f: \mathbb{N} \rightarrow \mathbb{R}^+$ be a function.
We define the $f$-\emph{kernel} to be
$$\Ker_f : = \{ i \in [n] \ : \ w_i \geq f(n) \}.$$
We will denote by $CL [\Ker_f]$ the subgraph of $CL (\bw)$ that is induced by the vertices of $\Ker_f$.
We show that there exists a function $f$ such that if we infect randomly $a(n)$ vertices  of $\Ker_f$,
then this is sufficient to infect almost the whole
of the $C$-kernel, for some constant $C>0$, with high probability. In other words, the gist of this theorem is that there is a specific part of the
random graph of size $o(n)$ such that if the initially infected vertices belong to it, then this
is enough to spread the infection to a positive fraction of the vertices.
\begin{theorem} \label{thm:fkernel}
Let $a: \mathbb{N} \rightarrow \mathbb{N}$ be a function such that $a(n) \rightarrow \infty$, as $n \rightarrow \infty$, but $a(n)=o(n)$.
Assume also ${r-1 \over 2r - \beta + 1} < \zeta \leq {1\over \beta -1}$.
If $\beta \in (2,3)$, then there exists an $\eps_0 = \eps_0 (\beta, \gamma_1, \gamma_2)$
such that for any positive $\eps < \eps_0$ there exists a constant $C=C(\gamma_1,\gamma_2, \beta, \eps, r) >0$ and
a function $f: \mathbb{N} \rightarrow \mathbb{R}^+$ such that $f(n) \rightarrow \infty$ as $n \rightarrow \infty$ but
$f(n) \ll n^{\zeta}$ satisfying the following. If we infect randomly $a(n)$ vertices in $\Ker_f$, then
at least $(1-\eps)|\Ker_C|$ vertices in $\Ker_C$ become infected a.a.s.
\end{theorem}
In both theorems, the sequence of probability spaces we consider are the product spaces of the random graph together with the
random choice of $\Ao$.

\section{Proofs}
In this section we present the proofs of Theorems~\ref{thm:global} and~\ref{thm:fkernel}.
We begin with stating a recent result due to Janson, {\L}uczak, Turova and Vallier~\cite{ar:JLTV10} regarding the evolution of bootstrap
percolation processes on Erd\H{o}s-R\'enyi random graphs, as these will be needed in our proofs.
These results regard the binomial model $G(N,p)$ introduced by Gilbert~\cite{ar:Gilb55} and subsequently became a major part of the
theory of random graphs
(see~\cite{Bol} or~\cite{JLR}).
Here $N$ is a natural number and $p$ is a real number that belongs to $[0,1]$. We consider the
set $[N]=: \{1,\ldots, N\}$ and create a random graph on the set $[N]$, including each pair $\{i,j\}$, where $i \not = j \in [N]$,
independently with probability $p$.

We begin with a few definitions as they were given in~\cite{ar:JLTV10}.
Recall that $r\geq 2$ is an integer and denotes the activation threshold.
We set
\begin{eqnarray}  \label{eq:ERFormulae}
&& T_c (N,p): = \left( {(r-1)!\over N p^r} \right)^{1/(r-1)}, \  A_c (N):= \left(1-{1\over r} \right)T_c(N,p), \ \\ \ \mbox{and} \ \nonumber
&& B_c (N):= N {(pN)^{r-1}\over (r-1)!} e^{-pN}.
\end{eqnarray}
Observe that if $p=p(N) \gg 1/N$, then $B_c(N) = o(N)$.
The following theorem is among the main results in~\cite{ar:JLTV10}.
\begin{theorem}[Theorem 3.1~\cite{ar:JLTV10}] \label{thm:JLTV}
Let $a : \mathbb{N} \rightarrow \mathbb{N}$ be a function. Assume that $\Ao$ is a subset of $[N]$ that
has size $a(N)$. Let $p=p(N)$ be such that $N^{-1} \ll p \ll N^{-1/r}$. Then a.a.s. (on the product space of the choice of $\Ao$ and the
random graph $G(N,p)$) we have
\begin{enumerate}
\item[(i)] if $a(N)/A_c(N) \rightarrow \alpha < 1$, then $|\Af| = (\phi (\alpha)+ o_p(1))T_c(N,p)$, where $\phi (\alpha)$ is the unique root in $[0,1]$ of
$$ r \phi (\alpha) - \phi(\alpha)^r  = (r-1) \alpha.$$
Further, $|\Af|/a(N) \stackrel{p}{\rightarrow} \phi_1 (\alpha):= {r \over r-1} \phi(\alpha)/\alpha$, with $\phi_1 (0):=1$;
\item[(ii)] if $a(N)/A_c(N) \geq 1 + \delta$, for some $\delta >0$, then $|\Af| = N - O_p (B_c(N))$. In other words, we have almost
complete percolation with high probability.
\item[(iii)] In case (ii), if further $a(N)\leq N/2$, we have complete percolation, that is, $|\Af|=N$ a.a.s., if and only if
$B_c (N) \rightarrow 0$, as $N\rightarrow \infty$, that is, if and only if $Np - (\log N + (r-1)\log \log N) \rightarrow \infty$ as
$N \rightarrow \infty$.
\end{enumerate}
\end{theorem}
The following theorem (also from~\cite{ar:JLTV10}) treats the dense regime.
\begin{theorem}[Theorem 5.8~\cite{ar:JLTV10}] \label{thm:JLTVDense}
Let $r\geq 2$. If $p\gg N^{-1/r}$ and $a(N)\geq r$, then a.a.s. $|\Af | = N$.
\end{theorem}
Now we proceed with the proof of Theorem~\ref{thm:fkernel}, as parts of it will be used in the proof of Theorem~\ref{thm:global}.

\subsection{Proof of Theorem~\ref{thm:fkernel}}
We will determine an $f$ sufficiently fast growing such that $CL[\Ker_f]$ stochastically ``contains" a dense enough $G(|\Ker_f|,p)$.
The density is high enough so that $a(n)$ exceeds the threshold given in Theorem~\ref{thm:JLTV} and therefore with
high probability we have the almost complete infection of $\Ker_f$. To show that the infection spreads over most of the vertices
of the $C$-kernel, we will split the set of vertices of $\Ker_C \setminus \Ker_f$ into ``bands"
$\La_j := \{ i \in [n] \ : \   f_j (n) \leq  w_i <  f_{j-1} (n) \}$, for $j=1,\ldots, T (n)$, where $T (n)$ as well as the functions
$f_j $ will be defined during our proof. We will show inductively that given that $\La_{j-1}$ is almost completely infected, then
with high probability $\La_j$ becomes almost completely infected as well. We now proceed with the details of the proof.

\subsubsection*{\bf Determining the function $f$}
We set
\begin{equation} \label{eq:fchoice}
f(n) = \left[ {(r-1)! W_{[n]}^r \over \gamma_1 n a^{r-1}(n)} \right]^{1 \over 2r - \beta +1}.
\end{equation}
We first need to show that $f(n) = o\left( n^{\zeta} \right)$, in order to ensure that $\Ker_f \not =\emptyset$.
\begin{claim} \label{clm:f1st}
If $\beta < 3$, then  $f(n) = o \left( n^{\zeta} \right)$.
\end{claim}
\begin{proof}
Note that
$$ f(n) = \Theta \left( \left({n \over a(n)} \right)^{r-1 \over 2r - \beta + 1}\right).$$
Thus the assumption that ${r-1 \over 2r - \beta + 1} < \zeta$ implies the claim.

\end{proof}
Since ${r-1 \over 2r -\beta +1} < {1\over 2}$, for $\beta < 3$, the proof of the above claim implies the following:
\begin{corollary}
If $\beta <3$, we have $f(n) = o (n^{1/2})$.
\end{corollary}
We also need to show that $a(n) \leq N_f$.
Recall that by Definition~\ref{def:Distr}, for all $n \in \mathbb{N}$ that are sufficiently large we have
\begin{equation} \label{eq:KerSize}
\gamma_1 f(n)^{-\beta + 1} \leq {N_f \over n} \leq \gamma_2 f(n)^{-\beta + 1}.
\end{equation}
\begin{claim} \label{clm:f2nd}
If $\beta < 3$, then $a(n) \ll N_f$.
\end{claim}
\begin{proof}
Observe that $a(n) = \Theta \left( n f^{{\beta -1 - 2r\over r-1}}(n) \right)$ whereas
$N_f = \Theta (n f^{-\beta + 1} (n))$. But
$$ {\beta - 1 - 2r \over r-1} < -\beta + 1,$$
which holds since $\beta < 3$.
This concludes the proof of the claim.
\end{proof}

For any two distinct vertices $i,j \in \Ker_f$ the probability of the edge
$\{ i,j\}$ being present is at least $p_f:=f^2 (n) / W_{[n]}$. Therefore setting $N_f : = |\Ker_f|$, it follows that $CL[\Ker_f]$
stochastically contains
$G(N_f,p_f)$. More precisely, there exists a coupling between these two random graphs such that always
$G (N_f, p_f) \subseteq CL [\Ker_f]$. Thus, if $\mathcal{P}$ is a non-decreasing property of graphs (that is, a set of graphs closed
under automorphisms and under the addition of edges), we have
\begin{equation} \label{eq:KernelDom}
\Pro { G(N_f, p_f) \in \mathcal{P}} \leq \Pro {CL [\Ker_f] \in \mathcal{P}}.
\end{equation}

To apply Theorem~\ref{thm:JLTV}, we first need to show that
$$ p_f \ll N_f^{-1/r}.$$
Let us set for convenience $x_{\beta,r}={r-1 \over 2r - \beta +1}$.
We have
$$p_f = \Theta \left( n^{2x_{\beta,r}-1} a^{-2x_{\beta,r}}(n) \right) $$
and $$N_f^{-1/r}:= \Theta \left( n^{-1/r+{(\beta -1)x_{\beta,r}\over r}} a^{-{(\beta -1)x_{\beta,r}\over r}} (n) \right).$$
In fact, we have
$$ 2x_{\beta,r}-1 = -1/r+{(\beta -1)x_{\beta,r}\over r}, $$
and
$$-2 x_{\beta,r} <  -{(\beta -1)x_{\beta,r}\over r},$$
as $2r \geq 4 > \beta -1$.

Now, from (\ref{eq:ERFormulae}) we have
$$ T_c (N_f, p_f) = \left( {(r-1)!\over N_f p_f^r} \right)^{1\over r-1}.$$
We will show that the choice of $f$ is such that $a(n) \geq T_c(N_f,p_f)$. Then since $A_c(N_f) = (1-1/r) T_c(N_f,p_f)$,
it follows from Theorem~\ref{thm:JLTV}(ii) together
with (\ref{eq:KernelDom}), that a.a.s. $CL [\Ker_f]$ becomes almost completely infected. To verify this statement,
we present the related calculations. By the choice of $f$ we have
$$a (n) = \left( {(r-1)! W_{[n]}^r \over n \gamma_1 f(n)^{-\beta +1+2r} } \right)^{1\over r-1} \geq T_c (N_f,p_f). $$

Finally, we need to bound $B_c(N_f)$ and show that $B_c(N_f) = o(N_f)$. To this end, it suffices to show that
$N_fp_f \rightarrow \infty$ as $n \rightarrow \infty$.
\begin{claim} \label{clm:f3rd}
If $\beta < 3$, then the function $f(n)$ is such that $N_f p_f \rightarrow \infty$ as $n \rightarrow \infty$.
\end{claim}
\begin{proof}
Indeed, by (\ref{eq:KerSize}) we have
\begin{equation*}
N_f p_f \geq n \gamma_1 {f(n)^{-\beta + 3} \over W_{[n]}} = \Omega \left( f(n)^{3-\beta}\right).
\end{equation*}
But $f(n) = \Theta \left( \left( {n\over a(n)}\right)^{r-1 \over 2r-\beta +1}\right)$ and, since $a(n) = o(n)$, this implies that
$f(n) \rightarrow \infty$ as $n \rightarrow \infty$.
This concludes the proof of the claim.
\end{proof}
By the above claims, Theorem~\ref{thm:JLTV}(ii) implies that for any $\eps >0$ a.a.s. at least $(1-\eps)|\Ker_f|$ vertices of $\Ker_f$
become infected.

\subsubsection*{\bf The dissemination of the infection in $\Ker_C \setminus \Ker_f$}
In this part of the proof we will show the following proposition which implies Theorem~\ref{thm:fkernel}.
\begin{proposition} \label{prop:InfectionSpread}
Let $f: \mathbb{N} \rightarrow \mathbb{R}^+$ be a function such that
$f(n) \rightarrow \infty$ as $n \rightarrow \infty$.
Then there exists an $\eps_0 = \eps_0 (\beta, \gamma_1, \gamma_2) >0$ such that for any positive $\eps < \eps_0$
there exists $C=C(\gamma_1,\gamma_2, \beta, \eps, r) >0$ for which the following holds.
If $(1-\eps)|\Ker_f|$ vertices of $\Ker_f$ have been infected, then a.a.s. at least $(1-\eps)|\Ker_C|$ vertices of
$\Ker_C$ become infected.
\end{proposition}
\begin{proof}
We will define a partition on the set of vertices in $\Ker_C \setminus \Ker_f$ as follows.
Firstly, we define a sequence of functions $f_i : \mathbb{N} \rightarrow \mathbb{R}^+$, for any integer $i \geq 0$.
Next, we define the real-valued function $\psi$ on the set of natural numbers
where for any $n \in \mathbb{N}$ we set $\psi(n) := {\ln C \over \ln f_0(n)}$ and  $C \in \mathbb{R}^+$ will be
determined during our proof.
Let $g_{\beta, n}(x) =(\beta -2)x + \psi(n)$, for $x \in \mathbb{R}$ and
for $i\in \mathbb{N}$ we let $g_{\beta, n}^{(i)}(x)$ be the $i$th iteration of $g_{\beta, n}(x)$ on itself.
We set $f_0 := f$ and for $i\geq 1$ we set $f_i:= f_{0}^{g_{\beta, n}^{(i)}(1)}$.
Also, for any $n \in \mathbb{N}$ we let
$T(n)$ be the maximum $i \in \mathbb{N}$ such that $f_0(n)^{g_{\beta, n}^{(i)}(1)} \geq C^{2\over 3-\beta}$.
Note that
\begin{equation} \label{eq:TBound} T(n) = O \left( \log \log n \right). \end{equation}

We are now ready to define the partition of $\Ker_C \setminus \Ker_f$.
For any $n \in \mathbb{N}$ and for $j = 1,\ldots, T(n)$ we set $$ \La_j : = \{i \in [n] \ : \  f_j(n) \leq w_i < f_{j-1} (n) \}.$$
We also set $\La_0 : = \Ker_f$ and finally $\La_{T(n)+1}:= \Ker_C \setminus \{ \cup_{j=0}^{T(n)} \La_j \}$.

The previous analysis has shown that a.a.s. if we infect randomly $a(n)$ vertices in $\La_0$, then at least $(1-\eps)|\La_0 |$ of $\La_0$ become
infected.
The following lemma serves as the inductive step, the proof of which is
postponed to the end of this section.
\begin{lemma} \label{lem:InductiveInfection}
There exists an $\eps_0 = \eps_0 (\beta, \gamma_1, \gamma_2) >0$ such that for any positive $\eps < \eps_0$
there exists $C=C(\gamma_1,\gamma_2, \beta, \eps, r) >0$ for which the following holds.
For $j=0,\ldots, T(n)$, if $(1-\eps)|\La_s|$ vertices of $\La_s$ have been infected, for $s=0,\ldots, j$, then with probability at least
$1 -  \exp \left(- {\eps^2 |\La_{j+1}|} \right)$ at least $(1- \eps) |\La_{j+1}|$ vertices of
$\La_{j+1}$ become infected.
\end{lemma}
The above lemma implies that the probability that for $j=1,\ldots, T(n)$ there are at least $(1-\eps)|\La_j|$ vertices in $\La_j$
that become infected, conditional on almost complete infection of $\Ker_f$, is at least
\begin{equation} \label{eq:ProbBound}
\begin{split}
\prod_{j=0}^{T(n)-1} \left( 1 -  \exp \left(- {\eps^2 |\La_{j+1}|\over 16} \right) \right)
\geq 1 - \sum_{j=0}^{T(n)-1} \exp \left(- {\eps^2 |\La_{j+1}| \over 16} \right).
\end{split}
\end{equation}
Thus, we need to bound from above the sum in the right-hand side of the above inequality.
To this end, we need to bound $|\La_{j+1}|$ from below, for $j=0,\ldots, T(n)-1$.
Note that
$$|\La_{j+1}| = |\Ker_{f_{j+1}}| - |\Ker_{f_{j}}|.$$
To bound the quantities on the right-hand side we use the bounds given by Definition~\ref{def:Distr}.
In particular, this implies that
$$|\La_{j+1}| \geq n \left( \gamma_1 f_{j+1}^{-\beta + 1}(n) - \gamma_2 f_{j}^{-\beta+1}(n)\right).$$
We use the definition of $f_j$ and in particular the identity:
\begin{eqnarray} \label{eq:frecursion}
f_{j+1}(n) &=& f_0^{g_{\beta,n}^{(j+1)}(1)}(n) = f_0^{(\beta -2)g_{\beta,n}^{(j)}(1) + \psi (n)}(n) \\
\nonumber &=& f_{j}^{\beta -2 }(n) f_0^{\psi (n)}(n) = f_{j}^{\beta -2}(n) C.
\end{eqnarray}
Thus we write
\begin{equation*}
\begin{split}
|\La_{j+1}| & \geq n \left( \gamma_1 f_{j+1}^{-\beta + 1}(n) - \gamma_2 f_{j}^{-\beta+1}(n)\right) \\
&= n \gamma_1 C^{-\beta + 1} f_{j}^{-(\beta - 1)(\beta -2)} (n)\left( 1- {\gamma_2 C^{\beta -1}\over \gamma_1}
f_{j}^{(\beta - 1)(\beta - 3)}(n)\right) \\
&\geq n \gamma_1 C^{-\beta + 1} f_{j}^{-(\beta - 1)(\beta -2)} (n)\left( 1- {\gamma_2 C^{-\beta +1}\over \gamma_1}
\right),
\end{split}
\end{equation*}
where in the last inequality we used the fact that $f_{j} (n) \geq C^{2 \over 3-\beta}$, which implies that
$ f_{j}^{(\beta -1)(\beta -3)} (n) \leq C^{-2(\beta -1)}$. Also, $f_{j} (n) \leq f_0(n) = f(n)$.
Thus, if $C$ is large enough, we obtain:
\begin{equation*}
\begin{split}
|\La_{j+1}| \geq  n f^{-(\beta - 1)(\beta -2)} (n) {\gamma_1 C^{-\beta + 1}  \over 2}.
\end{split}
\end{equation*}
But since $f(n)=o(n^{\zeta})$ we have  $f(n)= o(n^{1/(\beta -1)})$ as well.  So, for $n$ sufficiently large, $f(n) \leq n^{1/(\beta -1)}$.
Thereby we obtain $f^{(\beta-1)(\beta -2)}(n) \leq n^{\beta -2}$. Hence, there exists a constant $C'=C'(C,\gamma_1,\gamma_2)$ such that
for any $n$ sufficiently large we have
$$ |\La_{j+1}| \geq C' n^{3-\beta}.$$
Substituting this lower bound into the right-hand side of (\ref{eq:ProbBound}) and using (\ref{eq:TBound}) to bound the
number of summands there, we deduce that this sum is $o(1)$.

\end{proof}
\noindent
We conclude with the proof of Lemma~\ref{lem:InductiveInfection}.
\begin{proof}[Proof of Lemma~\ref{lem:InductiveInfection}]
Assume that for some integer $0\leq j \leq T(n)$, for $s=0,\ldots, j$ there are $(1-\eps)|\La_s|$ infected vertices in $\La_s$.
For $s=0,\ldots, j$ let $I_s \subset \La_s$ be the set of infected vertices of $\La_s$ and let $\hat{I}_{j} := \cup_{s=0}^{j} I_s$.

Consider a vertex $i \in \La_{j+1}$ and let $d_{\hat{I}_j}(i)$ denote the degree of $i$ in the set $\hat{I}_j$, that is, the number
of neighbours of vertex $i$ in this set.
We condition on the event that for $s=0,\ldots, j$ there are $(1-\eps)|\La_s|$ infected vertices in $\La_s$ -- thus every
probability calculation in this proof is conditional on this event.
We will first calculate the expected value of $d_{\hat{I}_j}(i)$ and show that it is large enough so
that the probability that $d_{\hat{I}_j}(i) < r$ is less than $\eps^2$. Thereafter, we apply the Chernoff bound for sums of indicator random
variables to bound the probability that there are at least $2\eps |\La_{j+1}|$ such vertices in the set $\La_{j+1}$ and conclude the proof of the
lemma.

To carry out these calculations we will need estimates on the total weight of the vertices that belong to a kernel.
Here and elsewhere the Landau notation involves absolute constants depending only on
$\gamma_1, \gamma_2$ as well as the average degree of the random graph.
\begin{claim} \label{lem:sumweights}
Let $f: \mathbb{N} \rightarrow \mathbb{R}^+$ be such that $f(n) \ll n^{\zeta}$.
For $\beta \in (2,3)$ we have
$$ \sum_{i \in \Ker_f } w_i = \Theta \left( {n \over f^{\beta -2}(n)}\right).$$
\end{claim}
\begin{proof}
By Definition~\ref{def:Distr}, there exists a positive real $x_0$ such that for every
$x_0 \leq s \leq n^{\zeta}$ we have
\begin{equation} \label{eq:defI} \gamma_1 s^{-\beta + 1} \leq 1 - F_n(s) \leq \gamma_2 s^{-\beta +1},\end{equation}
whereas for $s < x_0$ we have $F_n(s)=0$ and for $s > n^{\zeta}$ we have $F_n(s)=1$.
Thus, since $f(n) \leq n^{\zeta}$, we have
\begin{equation}\label{eq:KerSizeI} \gamma_1 f^{-\beta + 1}(n)\leq {|\Ker_f|\over n}  \leq \gamma_2 f^{-\beta +1}(n).
\end{equation}


We define the function $g_n$ on $[0,1]$ as follows. For $0 \leq x \leq 1- F_n (n^{\zeta})$ we set $g_n (x) = n^{\zeta}$ and
for $1-F_n(n^{\zeta}) < x \leq 1$ we set $g_n(x) = [1-F_n]^{-1}(x)$.
Thus we have
\begin{equation*}
\begin{split}
 \sum_{i \in \Ker_f} w_i & =n\int_{0}^{|\Ker_f|/n}g_n(x) dx = n \left( \int_{0}^{1-F_n(n^{\zeta})}g_n(x) dx +
\int_{1-F_n(n^{\zeta})}^{|\Ker_f|/n}g_n(x) dx\right) \\
& = n^{1+\zeta}(1-F_n(n^{\zeta})) + n \int_{1-F_n(n^{\zeta})}^{|\Ker_f|/n}g_n(x) dx \\
&= \Theta \left( n^{1-\zeta(\beta-2)} \right) + n \int_{1-F_n(n^{\zeta})}^{|\Ker_f|/n}g_n(x) dx .
\end{split}
\end{equation*}
Since $f(n) \ll n^{\zeta}$ it suffices to show that the second integral on the right-hand side satisfies the bounds of the claim.

Let us also define for every $x \in (0,1]$ the functions
$g_{1,n} (x) = \inf \{ s \ : \ \gamma_1 s^{-\beta + 1} \leq x \}$ and $g_{2,n} (x) = \inf \{ s \ : \ \gamma_2 s^{-\beta+1} \leq x \}$.
By (\ref{eq:defI}), for any $x \in (1-F_n(n^{\zeta}),1]$
$$\{ s \ : \ \gamma_2 s^{-\beta+1} \leq x \}
\subseteq \{ s \ : \ 1-F_n(s) \leq x \} \subseteq \{ s \ : \ \gamma_1 s^{-\beta + 1} \leq x \},$$
which implies that
$$ g_{1,n} (x) \leq g_n(x) \leq g_{2,n} (x). $$

Note that $g_{1,n} (x) = \left( {\gamma_1/x}\right)^{1\over \beta -1}$ and $g_{2,n} (x) = \left( {\gamma_2/x}\right)^{1\over \beta -1}$.
Hence
\begin{equation} \label{eq:sandwichI}
\int_{1-F_n(n^{\zeta})}^{|\Ker_f|/n} \left( {\gamma_1 \over x}\right)^{1\over \beta -1}dx  \leq
\int_{1-F_n(n^{\zeta})}^{|\Ker_f|/n} g_n(x) dx \leq
\int_{1-F_n(n^{\zeta})}^{|\Ker_f|/n} \left( {\gamma_2 \over x}\right)^{1\over \beta -1}dx.
\end{equation}
For $\ell\in \{1,2\}$ and since $\beta \in (2,3)$ we have
\begin{equation*}
\begin{split}
\int_{1-F_n(n^{\zeta})}^{|\Ker_f|/n} \left( {\gamma_\ell \over x}\right)^{1\over \beta -1}dx &= \gamma_\ell^{1\over \beta -1}
\int_{1-F_n(n^{\zeta})}^{|\Ker_f|/n} \left( {1\over x}\right)^{1\over \beta -1}dx \\
& = \gamma_\ell^{1\over \beta -1}~{\beta - 1 \over \beta-2}~ \left[\left( {|\Ker_f|\over n} \right)^{\beta -2 \over \beta -1} -
(1-F_n(n^{\zeta}))^{{\beta -2 \over \beta -1}}\right].
\end{split}
\end{equation*}
As $1- F_n (n^{\zeta}) = \Theta \left( n^{-\zeta(\beta -1)} \right)$ and $f(n)=o(n^{\zeta})$,
substituting the bounds of (\ref{eq:KerSizeI}) the claim follows.
\end{proof}

Let us continue with the estimate on $\Ex {d_{\hat{I}_j} (i)}$.
\begin{lemma} \label{lem:DegExpLow}
There exists $\eps_0 = \eps_0 (\beta, \gamma_1, \gamma_2)>0$ such that for every $\eps < \eps_0$ we have that
uniformly for all $i\in \La_{j+1}$ we have
\begin{equation*}
\Ex {d_{\hat{I}_j} (i)}  = \Theta (C).
\end{equation*}
\end{lemma}
\begin{proof}
Firstly, the definition of the random graph model implies that
\begin{equation} \label{eq:DegExpect}
\Ex {d_{\hat{I}_j} (i)} = {w_i W_{\hat{I}_j} \over W_{[n]}}.
\end{equation}
The weight $w_i$ is bounded from below by $f_{j+1}(n)$. We now need to bound from below $W_{\hat{I}_j}$.
\begin{claim}
We have
$$ W_{\hat{I}_j} = W_{\Ker_{f_j}} \left( 1 - O\left( \eps^{\beta- 2 \over \beta -1} \right) \right).$$
\end{claim}
\begin{proof}
Note that $|\hat{I}_j|\geq (1-\eps) |\Ker_{f_j}|$, which implies that $W_{\hat{I}_j}$ is at least as large
as the total weight of the $(1-\eps) |\Ker_{f_j}|$ vertices of smallest weight in $\Ker_{f_j}$.
To this end, we need to determine a function $\tilde{f}_j$ such that $\Ker_{\tilde{f}_j}$ has size
$\eps |\Ker_{f_j}|$. By Definition~\ref{def:Distr} we have
$$\gamma_1 \tilde{f}_j(n)^{-\beta + 1} \leq {|\Ker_{\tilde{f}_j}|\over n} = {\eps |\Ker_{f_j}|\over n}
\leq {\eps \gamma_2} f_j(n)^{-\beta +1}. $$
This inequality implies that
$$\tilde{f}_j(n) \geq f_j(n) \left( {1\over \eps}~{\gamma_1 \over \gamma_2}\right)^{1\over \beta -1}.$$
Thus, by Claim~\ref{lem:sumweights} we have
$$ W_{\Ker_{\tilde{f}_j}} = O \left( {n \over f_j^{\beta -2}(n)}
\left( {\eps \gamma_2 \over \gamma_1}\right)^{\beta-2 \over \beta -1}\right).$$
In turn, this implies that
$$ W_{\hat{I}_j} \geq W_{\Ker_{f_j}} \left( 1 - O\left( \eps^{\beta- 2 \over \beta -1} \right) \right).$$
As $W_{\hat{I}_j} \leq W_{\Ker_{f_j}}$, the claim follows.
\end{proof}
So if $\eps$ is small enough, then the right-hand side of (\ref{eq:DegExpect}) is bounded from below as follows.
\begin{equation}\label{eq:DegExpLow}
\begin{split}
\Ex {d_{\hat{I}_j} (i)}  \geq {1\over 2} {f_{j+1}(n)  W_{\Ker_{f_j}} \over W_{[n]}} = \Omega \left( f_{j+1}(n) f_j^{-\beta +2}(n)\right),
\end{split}
\end{equation}
using Claim~\ref{lem:sumweights} for the lower bound on $W_{\Ker_{f_j}}$.
By (\ref{eq:frecursion})
$$  f_{j+1}(n) f_j^{-\beta +2}(n) =f_j^{\beta -2 }(n) f_0^{\psi (n)}(n)  f_j^{-\beta +2}(n) = f_0^{\psi (n)}(n) = C.$$
Substituting this into the right-hand side of (\ref{eq:DegExpLow}) yields the lower bound in the lemma.
The upper bound also follows, observing that
$$ \Ex {d_{\hat{I}_j} (i)}  \leq {f_{j+1}(n)  W_{\Ker_{f_j}} \over W_{[n]}}. $$
\end{proof}
The next step is to show that if $C$ is large enough, then the probability that
$d_{\hat{I}_j}(i) < r$ can become as small as we need.
\begin{lemma} \label{lem:DegDistr}
Let $\eps_0$ be as in Lemma~\ref{lem:DegExpLow}.
For all $\eps < \eps_0$ there exists $C=C(\gamma_1, \gamma_2, \beta, \eps, r)>0$ such that for $j=0,\ldots, T(n)$ and for all $i \in \La_{j+1}$
$$\Pro {d_{\hat{I}_j}(i) < r } < \eps^2.$$
\end{lemma}
\begin{proof}
We will bound this probability using Chebyschev's inequality.
Observe that $d_{\hat{I}_j}(i) = \sum_{\ell \in \Ker_{f_j}} \Be \left( {w_i w_{\ell} \over W_n}\right)$, where the summands are
independent random variables. Hence, $\Var [d_{\hat{I}_j}(i) ] \leq \Ex {d_{\hat{I}_j}(i) }$.
By Lemma~\ref{lem:DegExpLow}, the latter is $\Theta (C)$. Thus, if $C$ is large enough, so that
$\Ex {d_{\hat{I}_j}(i)} -r > {C \over 2}$, then Chebyschev's inequality implies that
$$\Pro {d_{\hat{I}_j}(i) < r } = O \left( {1\over C} \right).$$
Making $C$ even larger, so that the right-hand side becomes smaller than $\eps^2$, the lemma follows.

\end{proof}
Setting $X_{j+1}:= \sum_{i \in \La_{j+1}} \mathbf{1}[d_{\hat{I}_j}(i) \geq r]$, we have
$I_{j+1}\geq X_{j+1}$.
As $X_{j+1}$ is the sum of independent identically distributed Bernoulli random variables we will use the Chernoff bound
to show that with high probability $X_{j+1} > |\La_{j+1}| (1-2\eps)$.
By Lemma~\ref{lem:DegDistr}
\begin{equation} \label{eq:InfectedExpectation}
\Ex {X_{j+1}} \geq |\La_{j+1}|(1-\eps^2).
\end{equation}
Thus, if $X_{j+1} \leq |\La_{j+1}| (1-2\eps)$, then $\Ex {X_{j+1}} - X_{j+1} \geq \eps |\La_{j+1}|$.
Setting $t:= \eps |\La_{j+1}|$, Theorem~2.1 in~\cite{JLR} (Inequality (2.6))  yields
\begin{equation*}
\begin{split}
\Pro {X_{j+1} \leq \Ex {X_{j+1}} - t} \leq \exp \left( - {t^2 \over 2 \Ex {X_{j+1}}} \right) \leq \exp
\left(- {\eps^2\over 2} |\La_{j+1}| \right),
\end{split}
\end{equation*}
for $\eps < 1/2$.
Thus, we deduce that
\begin{equation} \label{eq:azumafinal}
\begin{split}
\Pro {I_{j+1}\geq |\La_{j+1}| (1-2\eps)} \geq 1 - \exp \left(- {\eps^2\over 2} |\La_{j+1}| \right).
\end{split}
\end{equation}

\end{proof}

\subsection{Proof of Theorem~\ref{thm:global}: Subcritical case}
We will use a first moment argument to show that if $a(n) =o (a_c(n))$, then a.a.s. there are no
vertices outside $\Ao$ that have at least $r$ neighbours in $\Ao$ and, therefore, the bootstrap
percolation process does not actually evolve. Here we assume that initially each vertex becomes infected with
probability $a(n)/n$, independently of every other vertex.

For every vertex $i \in [n]$, we define an indicator random variable $X_i$ which is 1 precisely when
vertex $i$ has at least $r$ neighbours in $\Ao$. Let $X = \sum_{i \in [n]} X_i$. Our aim is to show that
$\Ex {X} = o(1)$, thus implying that a.a.s. $X=0$.

For $i\in [n]$ let $p_i=\Ex{X_i} = \Pro {X_i=1}$. We will first give an upper bound on $p_i$ and, thereafter, the
linearity of the expected value will conclude our statement.

\begin{lemma} \label{lem:pibound}
For all integers $r\geq 2$ and all $i\in [n]$, we have
$$p_i  \leq \left( {e w_i a(n) \over r n}\right)^r. $$
\end{lemma}
From this, we can use the linearity of the expected value to deduce an upper bound on $\Ex {X}$.
We have
\begin{equation} \label{eq:XExp}
\begin{split}
\Ex {X} = \sum_{i \in [n]} p_i \leq \sum_{i \in [n]} \left( {e w_i a(n) \over r n} \right)^r = o
 \left( \left({a_c (n) \over n}\right)^r \right)~\sum_{i \in [n]} w_i^r.
\end{split}
\end{equation}
We now need to give an estimate on $\sum_{i \in [n]} w_i^r$.
\begin{claim}
For all integers $r\geq 2$ and for $\beta \in (2,3)$ we have
$$ \sum_{i \in [n]} w_i^r = \Theta \left( n^{1+\zeta (r -\beta +1)} \right).$$
\end{claim}
\begin{proof}
By Definition~\ref{def:Distr}, there exists a positive real $x_0$ such that for every
$x_0 \leq s \leq n^{\zeta}$ we have
\begin{equation}\label{eq:def} \gamma_1 s^{-\beta + 1} \leq 1 - F_n(s) \leq \gamma_2 s^{-\beta +1},\end{equation}
whereas for $s < x_0$ we have $F_n(s)=0$ and for $s > n^{\zeta}$ we have $F_n(s)=1$.
As before, we define the function $g_n$ on $[0,1]$ as follows. For $0 \leq x \leq 1- F_n (n^\zeta)$ we set $g_n (x) = n^{\zeta}$ and
for $1-F_n(n^{\zeta}) < x \leq 1$ we set $g_n(x) = [1-F_n]^{-1}(x)$.
Hence, we write
\begin{equation*}
\begin{split}
 \sum_{i \in [n]} w_i^r & =n\int_{0}^{1}g_n^r(x) dx = n \left( \int_{0}^{1-F_n(n^{\zeta})} g_n^r (x) dx  +
\int_{1-F_n(n^{\zeta})}^1 g_n^r (x) dx \right) \\
& = \Theta \left( n^{1 + \zeta (r -\beta + 1)} \right) + n\int_{1-F_n(n^{\zeta})}^1 g_n^r (x) dx.
\end{split}
\end{equation*}
Hence, it suffices to show that the integral on the right-hand side satisfies the bounds of the claim.

Let us also define for every $x \in (0,1]$ the functions
$g_{1,n} (x) = \inf \{ s \ : \ \gamma_1 s^{-\beta + 1} \leq x \}$ and $g_{2,n} (x) = \inf \{ s \ : \ \gamma_2 s^{-\beta+1} \leq x \}$.
By (\ref{eq:def}), for any $x \in (1-F_n(n^{\zeta}),1]$
$$\{ s \ : \ \gamma_2 s^{-\beta+1} \leq x \}
\subseteq \{ s \ : \ 1-F_n(s) \leq x \} \subseteq \{ s \ : \ \gamma_1 s^{-\beta + 1} \leq x \},$$
which implies that
$$ g_{1,n} (x) \leq g_n(x) \leq g_{2,n} (x). $$

Note that $g_{1,n} (x) = \left( {\gamma_1/x}\right)^{1\over \beta -1}$ and $g_{2,n} (x) = \left( {\gamma_2/x}\right)^{1\over \beta -1}$.
Hence
\begin{equation} \label{eq:sandwich}
\int_{1-F_n(n^{\zeta})}^{1} \left( {\gamma_1 \over x}\right)^{r\over \beta -1}dx  \leq
\int_{1-F_n(n^{\zeta})}^{1} g_n^r(x) dx \leq
\int_{1-F_n(n^{\zeta})}^{1} \left( {\gamma_2 \over x}\right)^{r\over \beta -1}dx.
\end{equation}
For $\ell\in \{1,2\}$, since $\beta \in (2,3)$ and $r\geq 2$, we have
\begin{eqnarray*}
\int_{1-F_n(n^{\zeta})}^{1} \left( {\gamma_\ell \over x}\right)^{r\over \beta -1}dx &=& \gamma_\ell^{r\over \beta -1}
\int_{1-F_n(n^{\zeta})}^{1} \left( {1\over x}\right)^{r\over \beta -1}dx \\
&=& \gamma_\ell^{r\over \beta -1}~{\beta - 1 \over r-\beta+1}~\left[ (1-F_n(n^{\zeta}))^{-{r \over \beta -1}+1}
- 1 \right].
\end{eqnarray*}
Recall that $1-F_n(n^{\zeta}) = \Theta (n^{-\zeta(\beta -1)})$.
Thus through (\ref{eq:sandwich}) we deduce that for $r\geq 2$ and $\beta \in (2,3)$
$$ n \int_{1- F_n(n^{\zeta})}^{1} g^r(x) dx = \Theta \left( n^{1+\zeta(r - \beta + 1)}\right). $$
The claim now follows.
\end{proof}
Substituting this bound into the right-hand side of (\ref{eq:XExp}), we obtain:
\begin{equation*}
\begin{split}
\Ex {X} = o \left(  {n^{r(1-\zeta) + \zeta (\beta - 1) - 1} \over n^r}~n^{1+\zeta (r - \beta + 1)}\right).
\end{split}
\end{equation*}
But
$$ r(1-\zeta) + \zeta (\beta - 1) - 1 - r +1 +\zeta (r - \beta + 1)= 0, $$
thus implying that  $\Ex {X} = o(1)$.

We finish the proof of this part of Theorem~\ref{thm:global} with the proof of Lemma~\ref{lem:pibound}.
\begin{proof}[Proof of Lemma~\ref{lem:pibound}]
Note that for all $i\in [n]$ we have
$$ p_i = \Pro {\sum_{j \in [n] \setminus \{i\} } e_{ij} {\mathbf 1}[j \in \Ao] \geq r},$$
where $e_{ij}$ is the indicator random variable that is equal to 1 precisely when the pair $\{i,j\}$
belongs to the set of edges of $CL (\bw)$.
The random variable $e_{ij}{\mathbf 1}_{j \in \Ao}$ is  Bernoulli distributed with expected value
equal to ${w_i w_j \over W_{[n]}}~{a(n) \over n}$. We denote it by $I_j$, for all $j\in [n]\setminus \{ i\}$.

We will use a Chernoff-bound-like technique to bound this probability. Let $\theta >0$ be a real number.
We have
\begin{equation*}
\begin{split}
&\Pro {\sum_{j \in [n]\setminus\{ i \} } I_j \geq r} =
\Pro { \exp \left( \theta  \sum_{j \in [n]\setminus\{ i \} } I_j\right) \geq \exp \left( \theta r \right) } \\
& \leq {\Ex {  \exp \left( \theta\sum_{j \in [n]\setminus\{ i \} } I_j\right)} \over e^{\theta r}} = {\prod_{j \in [n]\setminus\{ i \}} \Ex {e^{\theta I_j}}\over e^{\theta r}}\\
& = {\prod_{j \in [n]\setminus\{ i \}} \left(e^{\theta}~{w_i w_j \over W_{[n]}}~{a(n) \over n} +
\left(1- {w_i w_j \over W_{[n]}}~{a(n) \over n} \right)\right) \over e^{\theta r}} \\
&\leq {\prod_{j \in [n]\setminus\{ i \}} \exp \left((e^{\theta} -1)~{w_i w_j \over W_{[n]}}~{a(n) \over n} \right)\over e^{\theta r}}
= {\exp \left((e^{\theta} -1)~\sum_{j \in [n]\setminus\{ i \}}{w_i w_j \over W_{[n]}}~{a(n) \over n} \right)\over e^{\theta r}} \\
&\leq {\exp \left((e^{\theta} -1)~ w_i~{a(n) \over n} -\theta r \right) }.
\end{split}
\end{equation*}
The exponent in the last expression is minimized when $\theta$ is such that $e^{\theta} = {r n \over w_i a(n)}$.
Thus, we obtain
\begin{equation*}
\begin{split}
\Pro {\sum_{j \in [n]\setminus\{ i \} } I_j \geq r} & \leq \exp \left( r - {w_i a(n) \over n}\right)~ \left( {w_i a(n)\over rn} \right)^r \\
& = \left[\exp \left( 1 - {w_i a(n) \over rn}\right)~ \left( {w_i a(n)\over rn} \right) \right]^r \leq \left( {e w_i a(n)\over rn} \right)^r.
\end{split}
\end{equation*}
\end{proof}

\subsection{Proof of Theorem~\ref{thm:global}: Supercritical case}
In this part of the proof, we shall be assuming that $a_c(n) = o(a(n))$.
Additionally, we shall assume that the initially infected set is the set of the $a(n)$ vertices
of smallest weight.

We will show first that there exists a function $f:\mathbb{N} \rightarrow \mathbb{R}^+$ such that
$f(n) \rightarrow \infty$ as $n \rightarrow \infty$ but $f(n) = o (n^{\zeta})$ for which a.a.s.
$\Ker_f$ will become completely infected.
Thereafter, using the proof of Theorem~\ref{thm:fkernel} we will
deduce that there exists a real number $C >0$ such that with high probability $\Ker_C$ will be almost completely infected.
This implies that there exists an $\eps >0$ such that a.a.s. at least $\eps n$ vertices become infected.

\subsubsection*{\bf Spreading the infection to a positive fraction of the vertices}

We begin with determining the function $f$ much as we did in the proof of Theorem~\ref{thm:fkernel}.
To this end, we need to bound from below the probability that an arbitrary vertex in $\Ker_f$ becomes
infected. In fact, we shall bound from below the probability that an arbitrary vertex in $\Ker_f$ will become
infected already in the first round. Note that this amounts to bounding the probability that such a vertex has at least $r$
neighbours in $\Ao$. Therefore, this forms a collection of independent events which is equivalent to the random independent
infection of the vertices of $\Ker_f$ with probability equal to the derived lower bound. Recall that the random graph induced
on $\Ker_f$ stochastically contains an Erd\H{o}s-R\'enyi random graph with the appropriate parameters. This observation allows
us to determine $f$. To be more specific, if the probability that any given vertex in $\Ker_f$ exceeds the complete infection threshold
of this  Erd\H{o}s-R\'enyi random graph and the condition of Theorem~\ref{thm:JLTV}(iii) is satisfied, then a.a.s. $\Ker_f$ eventually becomes completely infected. This condition will specify $f$.

Under the assumption that $\Ao$ consists of the $a(n)$ vertices of smallest weight, we will bound from below the probability
a vertex $v \in \Ker_f$ has at least $r$ neighbours in $\Ao$. We denote the degree of $v$ in $\Ao$ by $d_{\Ao}(v)$ and
note that this random variable is equal to $\sum_{i \in \Ao} \Be \left( {w_v w_i \over W_{[n]}} \right)$, where the summands
are independent Bernoulli distributed random variables.
Note also that for all $n$ and for all $i \in [n]$ we have $w_i \geq x_0$.
Thus, we can deduce the following (parts of it hold for $n$ sufficiently large)
\begin{equation*}
\begin{split}
\Pro {\sum_{i \in \Ao} \Be \left( {w_v w_i \over W_{[n]}} \right) \geq r} & \geq
\Pro {\sum_{i \in \Ao} \Be \left( {w_v x_0 \over W_{[n]}} \right) \geq r} \\ & =
\Pro { \Bin \left( a(n), {w_v x_0 \over W_{[n]}} \right) \geq r } \\ & \geq
{a(n) \choose r}~\left({w_v x_0 \over W_{[n]}} \right)^r ~\left( 1- {w_v x_0 \over W_{[n]}} \right)^{a(n)-r} \\ & \geq
{a(n)^r \over 1.5~r!}~\left({f(n) x_0 \over W_{[n]}} \right)^r ~\left( 1- {f(n)  x_0 \over W_{[n]}} \right)^{a(n)-r}.
\end{split}
\end{equation*}
Thus, assuming that $a(n)f(n) = o(n)$ we have
$$\left( 1- {f(n)  x_0 \over W_{[n]}} \right)^{a(n)-r} = 1-o(1).$$
Therefore, for $n$ sufficiently large
\begin{equation} \label{eq:InfectioPro}
\Pro {\sum_{i \in \Ao} \Be \left( {w_v w_i \over W_{[n]}} \right) \geq r} \geq {1\over 2 r!}~ \left({a(n) f(n) x_0 \over W_{[n]}} \right)^r
=:p_{Inf}.
\end{equation}
Thus every vertex of $\Ker_f$ becomes infected during the first round with probability at least $p_{Inf}$ independently of every
other vertex in $\Ker_f$.

We shall first consider the case where ${2r - \beta +1 \over r-1} \leq \zeta \leq {1\over \beta -1}$, where
$a_c (n) = n^{r(1-\zeta) + \zeta (\beta -1) - 1 \over r}$.
Let us assume that $a(n)=\omega (n) n^{r(1-\zeta) + \zeta (\beta -1) - 1 \over r}$, where
$\omega: \mathbb{N} \rightarrow \mathbb{R}^+$ is some increasing function that grows slower than any polynomial.
Setting $f=f(n)={n^{\zeta} \over \omega^{1+1/r} (n)}$, we will consider $CL[\Ker_f]$.
Before doing so, we will verify the assumption that $a(n)f(n)=o(n)$. Indeed, we have
\begin{equation*}
\begin{split}
a(n)f(n)= {1\over \omega^{1/r}(n)}~n^{{r(1-\zeta) + \zeta (\beta -1) - 1 \over r}+\zeta}.
\end{split}
\end{equation*}
But
\begin{equation*}
\begin{split}
&{r(1-\zeta) + \zeta (\beta -1) - 1 \over r}+\zeta = {r(1-\zeta) + \zeta (\beta -1) - 1 +r\zeta \over r}\\
 & =
1+ {\zeta (\beta - 1) - 1\over r} \leq 1,
\end{split}
\end{equation*}
since $\zeta \leq 1/(\beta -1)$, whereby $a(n)f(n) \leq {n \over \omega^{1/r}(n)} = o(n)$.

Now, note that if $\zeta > {1 \over 2}$, then $CL[\Ker_f]$ is the complete graph on $|\Ker_f|$ vertices.
However, when $\zeta \leq {1 \over 2}$, then  $CL [\Ker_f]$ stochastically contains $G(N_f,p_f)$, where $N_f = |\Ker_f|$ and
$p_f = {f^2(n) \over W_{[n]}}$. We will treat these two cases separately.

\medskip
\noindent
\emph{Case I}: ${1\over 2} < \zeta \leq {1\over \beta -1}$.

In this case, as $CL[\Ker_f]$ is the complete graph, it suffices to show that with high probability at least $r$ vertices
of $\Ker_f$ become infected already at the first round. In fact, we will show that the expected number of vertices of
$\Ker_f$ that become infected during the first round tends to infinity as $n$ grows.
Note that this number is equal to $N_f p_{Inf}$. Thus, once we show that $N_f p_{Inf} \rightarrow \infty$, as $n \rightarrow \infty$,
then Chebyschev's inequality or a standard Chernoff bound can show that with probability $1-o(1)$, there are at least
$r$ infected vertices in $\Ker_f$ and, thereafter, the whole of $\Ker_f$ becomes infected in one round.

We have
$$ N_f = |\Ker_f| = \Omega \left ( n \left( \omega(n) \over n^\zeta \right)^{\beta - 1}\right),$$
and by (\ref{eq:InfectioPro}) we have
$$ p_{Inf} = \Theta \left( {1\over \omega (n)}~\left( { n^{r(1-\zeta) + \zeta (\beta -1) -1\over r} \cdot n^{\zeta}\over n } \right)^r \right) =
\Theta \left( {n^{\zeta (\beta -1) - 1}\over \omega (n)} \right).
$$
Hence
$$ N_f p_{Inf} = \Omega \left(\omega^{\beta-2}(n)\right). $$

\medskip

\noindent
\emph{Case II}: ${r-1\over 2r - \beta + 1} < \zeta \leq {1 \over 2}$.

As we mentioned above, $CL [\Ker_f]$ stochastically contains $G(N_f,p_f)$, where $p_f = {f^2(n) \over W_{[n]}}$,
as $\zeta \leq {1\over 2}$. We will show that here $N_fp_f^r \rightarrow \infty$ as $n \rightarrow \infty$ and
by Theorem~\ref{thm:JLTVDense} we deduce that $\Ker_f$ becomes completely infected with probability $1-o(1)$.
We have
\begin{equation} \label{eq:T_c}
\begin{split}
N_fp_f^r = \Theta \left( \omega^{\beta-1}(n) n^{1-\zeta(\beta -1)}{n^{2\zeta r}\over \omega^{2r+2}(n) n^r} \right).
\end{split}
\end{equation}
 and the expression on the right-hand side is
 $$\omega^{-(2r-\beta+3)}(n) n^{-(r-1) + \zeta (2r-\beta +1)} \rightarrow \infty, $$
 by our assumption on $\zeta$.

\medskip
Finally, we deal with smaller values of $\zeta$, proving the last part of Theorem~\ref{thm:global}.

\medskip
\noindent
\emph{Case III}: $0< \zeta \leq {r-1 \over 2r -\beta +1}$.

In this case, we appeal to Theorem~\ref{thm:JLTV}.  We will show first that
$$ p_f \ll N^{-1/r}_f. $$
By (\ref{eq:T_c}) we have
$$ N_f p_f^r = \Theta \left( \omega^{-(2r- \beta +3)}(n) n^{1-\zeta(\beta -1) +r (2\zeta - 1)} \right). $$
The second exponent on the right-hand side of the above is equal to $- (r-1) + \zeta (2r -\beta +1) \leq 0$, by our assumption
on $\zeta$, whereby we have $p_f = o\left(N_f^{-1/r}\right)$.

Recall that  $a_c^+ (n) = n^{1 - \zeta {r-\beta +2 \over r-1}}$.
Let us set $\xi = 1 - \zeta {r-\beta +2 \over r-1}$.
It suffices to show that
\begin{equation} \label{eq:CriticalCond} N_f p_{Inf} \gg T_c (N_f, p_f).
\end{equation}
Since $T_c (N_f,p_f) = \Theta \left( \left( N_f p_f^r \right)^{-{1\over r-1}} \right)$, the above calculation implies that
$$ T_c (N_f,p_f) = \Theta \left( \omega^{2r -\beta + 3 \over r-1}(n) n^{1 - \zeta {2r-\beta + 1\over r-1}}\right).$$
Let $a(n) = \omega^2 (n) a_c^+ (n)$.
Then
$$ N_f p_{Inf} = \Theta \left( a^r(n) n^{1 - \zeta (\beta -1) + \zeta r - r }\right). $$
Hence
$$ N_f p_{Inf} = \Theta \left( \omega^{2r} (n) n^{r\xi -(r-1) + \zeta (r - \beta + 1)} \right).$$
But $\xi$ satsfies
$$ r \xi = r -\zeta \left( r-\beta + 1 +{2r - \beta +1 \over r-1} \right),$$
since
\begin{equation*}
\begin{split}
& r-\beta + 1 +{2r - \beta +1 \over r-1}  = {(r-1)(r-\beta + 1) + 2r -\beta + 1\over r-1} \\
& =  {r(r-\beta + 1) + r \over r-1} = r ~{r -\beta +2 \over r-1}.
\end{split}
\end{equation*}
Hence
$$ r \xi  -(r-1) + \zeta (r - \beta + 1) =  1- \zeta {2r -\beta + 1 \over r-1}. $$
Also ${2r -\beta + 3 \over r-1} \leq 2r - \beta + 3 < 2r$, since $r\geq 2$. Thus (\ref{eq:CriticalCond}) follows.

For each one of the above cases, Proposition~\ref{prop:InfectionSpread} implies that for any real $\eps >0$ that is small enough there exists a real number
$C=C(\gamma_1,\gamma_2,\beta, \eps) >0$ such that a.a.s. at least $(1-\eps)|\Ker_C|$ vertices
of $\Ker_C$ become infected. But by (\ref{eq:KerSize}) we have $|\Ker_C| = \Theta (n)$ and the second part of Theorem~\ref{thm:global}
follows.

\medskip
\noindent
{\bf Acknowledgment} We would like to thank Rob Morris for pointing an oversight in an earlier version of this paper.

\bibliographystyle{plain}
\bibliography{btsrp}

\end{document}